\documentclass{amsart}%
\usepackage{amssymb}
\usepackage{amsfonts}
\usepackage{amsmath}
\usepackage{graphicx}%
\providecommand{\U}[1]{\protect\rule{.1in}{.1in}}
\newtheorem{theorem}{Theorem}[section]
\theoremstyle{plain}

\newtheorem{corollary}[theorem]{Corollary}

\newtheorem{definition}[theorem]{Definition}
\newtheorem{example}[theorem]{Example}

\newtheorem{lemma}[theorem]{Lemma}

\newtheorem{proposition}[theorem]{Proposition}
\newtheorem{remark}[theorem]{Remark}

\numberwithin{equation}{section}

\begin{document}
\title[Fuchsian Groups actions on noncommutative Hardy algebras]{Ergodic Actions of Convergent Fuchsian groups on quotients of the
noncommutative Hardy algebras}
\author{Alvaro Arias}
\address{Department of Mathematics\\
University of Denver\\
Denver CO 80208}
\email{aarias@math.du.edu}
\urladdr{http://www.math.du.edu/\symbol{126}aarias}
\author{Fr\'{e}d\'{e}ric Latr\'{e}moli\`{e}re}
\email{frederic@math.du.edu}
\urladdr{http://www.math.du.edu/\symbol{126}frederic}
\date{February 10, 2010}

\begin{abstract}
We establish that particular quotients of the non-commutative Hardy algebras
carry ergodic actions of convergent discrete subgroups of the group
$\operatorname*{SU}(n,1)$ of automorphisms of the unit ball in $\mathbb{C}%
^{n}$. To do so, we provide a mean to compute the spectra of quotients of
noncommutative Hardy algebra and characterize their automorphisms in term of
biholomorphic maps of the unit ball in $\mathbb{C}^{n}.$

\end{abstract}
\maketitle

\bigskip We establish that given any discrete subgroup $\Gamma$ of
$\operatorname*{SU}(n,1)$ such that the orbit of $0$ for the action of
$\Gamma$ on the open unit ball $\mathbb{B}_{n}$ of $\mathbb{C}^{n}$ satisfies
the Blaschke condition:%
\[
\sum_{\gamma\in\Gamma}\left(  1-\left\Vert \gamma(0)\right\Vert _{\mathbb{C}%
^{n}}\right)  <\infty\text{,}%
\]
there exists a quotient algebra of the noncommutative Hardy algebra
$F_{n}^{\infty}$ whose group of weak* continuous automorphisms is the
stabilizer $\underline{\Gamma}$ of the orbit of $0$ for $\Gamma$ in
$\mathbb{B}_{n}$. Moreover, $\Gamma$ acts ergodically on this quotient
algebra. Our methods rely heavily on the theory of analytic functions in
several variables.

\bigskip The noncommutative Hardy algebra $F_{n}^{\infty}$ is the
weak-operator closure of the left regular representation of the free semigroup
on $n$ generators, and it is a noncommutative analogue of the Hardy algebra
$H^{\infty}\left(  \mathbb{B}_{1}\right)  $ \cite{Popescu91}. Our motivation
for this study is to explore the very rich ideal structure of $F_{n}^{\infty}%
$, as well as the structure of automorphism groups of non-self-adjoint
operator algebras. This paper is based upon the remarkable result of Davidson
and Pitts \cite[Theorem 4.11]{Pitts98}\ that the group of completely
contractive automorphisms of $F_{n}^{\infty}$ is $\operatorname*{SU}(n,1)$,
i.e. the same group as the group of biholomorphic maps from $\mathbb{B}_{n}$
into itself \cite[Lemma 4.9]{Pitts98}. See also \cite{Voiculescu85} and
\cite{Popescu10} for other approaches to this fact. We recall that
$\operatorname*{SU}(n,1)$ is the Lie group of $\left(  n+1\right)
\times\left(  n+1\right)  $ matrices of bilinear maps on $\mathbb{C}^{n}$
preserving the canonical sesquilinear form of signature $\left(  1,n\right)  $.

\bigskip Let us recall the definition of $F_{n}^{\infty}$ \cite{Popescu91}.
Let $n\in\mathbb{N}$ with $n>0$. The full Fock space $F^{2}\left(
\mathbb{C}^{n}\right)  $ is the completion of:%
\[%
{\displaystyle\bigoplus\limits_{k\in\mathbb{N}}}
\left(  \mathbb{C}^{n}\right)  ^{\otimes k}=\mathbb{C}\oplus\mathbb{C}%
^{n}\oplus\left(  \mathbb{C}^{n}\otimes\mathbb{C}^{n}\right)  \oplus\left(
\mathbb{C}^{n}\otimes\mathbb{C}^{n}\otimes\mathbb{C}^{n}\right)  \oplus\cdots
\]
for the Hilbert norm associated to the inner product $\left\langle
.,.\right\rangle $ defined on elementary tensors by:%
\[
\left\langle \xi_{0}\otimes\cdots\otimes\xi_{m},\zeta_{0}\otimes\ldots
\otimes\zeta_{k}\right\rangle =\left\{
\begin{array}
[c]{c}%
0\text{ if }m\not =k\text{,}\\%
{\displaystyle\prod\limits_{j=0}^{n}}
\left\langle \xi_{j},\zeta_{j}\right\rangle _{\mathbb{C}^{n}}\text{
otherwise,}%
\end{array}
\right.
\]
where $\left\langle .,.\right\rangle _{\mathbb{C}^{n}}$ is the canonical inner
product on $\mathbb{C}^{n}$. Let $\left\{  e_{1},\ldots,e_{n}\right\}  $ be
the canonical basis of $\mathbb{C}^{n}$ and, naturally, let $1$ be the
canonical basis of $\mathbb{C}$. We define, for each $j\in\left\{
1,\ldots,n\right\}  $, an operator $S_{j}$ on $F_{n}^{2}\left(  \mathbb{C}%
^{n}\right)  $ by:%
\[
S_{j}\left(  e_{i_{1}}\otimes\cdots\otimes e_{i_{m}}\right)  =e_{j}\otimes
e_{i_{1}}\otimes\cdots\otimes e_{i_{m}}\text{ and }S_{j}\left(  1\right)
=e_{j}\text{.}%
\]
The operators $S_{1},\ldots,S_{n}$ are called the left creation operators, and
we observe that $\sum_{j=1}^{n}S_{j}S_{j}^{\ast}\leq1$, i.e. $\left[
S_{1}\ \cdots\ S_{n}\right]  $ is a row contraction. The
weak-operator-topology closure of the algebra generated by $\left\{
1,S_{1},\ldots,S_{n}\right\}  $ is the noncommutative Hardy algebra
$F_{n}^{\infty}$. The fundamental property of $F_{n}^{\infty}$ is that, given
any separable Hilbert space $\mathcal{H}$ and any $n$-tuple $T=\left(
T_{1},\ldots,T_{n}\right)  $ of operators on $\mathcal{H}$ such that
$\sum_{j=1}^{n}T_{j}T_{j}^{\ast}<1$, there exists a unique completely
contractive algebra homomorphism $\pi_{T}$ from $F_{n}^{\infty}$ into the
algebra of bounded linear operators on $\mathcal{H}$ such that $\pi_{T}%
(S_{j})=T_{j}$ for $j=1,\ldots,n$, and we note that this map is weak*
continuous. This property was established by Popescu \cite{Popescu91} and in
\cite{Popescu99}\ using a noncommutative generalization of the Poisson
transform. As a matter of notation, we will write $\varphi\left(  T_{1}%
,\ldots,T_{n}\right)  $ for the operator $\pi_{T}\left(  \varphi\right)  $
whenever $\varphi\in F_{n}^{\infty}$. The algebra $F_{n}^{\infty}$ plays a
very important role in interpolation theory, among other matters, and many
results valid for the Hardy algebra $H^{\infty}\left(  \mathbb{B}_{1}\right)
$ can be extended to $F_{n}^{\infty}$.

\bigskip The Banach space $F_{n}^{\infty}$ is a dual space and thus can be
endowed with the corresponding weak*\ topology, which agrees with the
restriction of the weak-operator topology to $F_{n}^{\infty}$. In this paper,
we will refer to this topology on $F_{n}^{\infty}$ as the weak* topology on
$F_{n}^{\infty}$ \cite{Davidson99}. We also note that Davidson and Pitts use
the notation $\mathcal{L}_{n}$ for $F_{n}^{\infty}$ but we shall prefer
Popescu's notation to emphasize the connection with the Hardy algebra.

In this paper, we define the spectrum of an operator algebra as the space of
all weak* continuous scalar-valued algebra homomorphisms, endowed with the
$\sigma\left(  F_{n}^{\infty\ast},F_{n}^{\infty}\right)  $ topology (i.e. the
weak* topology on the dual of $F_{n}^{\infty}$). As a consequence of the
fundamental property of $F_{n}^{\infty}$, we note that the spectrum of
$F_{n}^{\infty}$ consists exactly of the maps $\pi:F_{n}^{\infty
}\longrightarrow\mathbb{C}$ such that $\sum_{j=1}^{n}\left\vert \pi
(S_{j})\right\vert ^{2}<1$, and one checks that indeed the spectrum is
homeomorphic to $\mathbb{B}_{n}$ \cite{Popescu96}.

While we use Popescu's $F_{n}^{\infty}$ in this paper, one should observe that
all our quotient algebras are in fact commutative. Indeed, we quotient
$F_{n}^{\infty}$ by intersections of kernels of scalar-valued algebra
homomorphisms, which always contain the commutator ideal of $F_{n}^{\infty}$.
Hence, all our constructions factor through the multiplier algebra of the
symmetric Fock space, as studied in \cite{Arveson98}.

\bigskip This paper is organized as follows. The first section provides a
functional mean to compute the pseudohyperbolic metric on $\mathbb{B}_{n}$
which will prove important for our purpose. We then compute the spectrum of a
large class of quotient algebras of $F_{n}^{\infty}$. We then conclude that
some of these quotients admit a discrete subgroup of $SU(n,1)$ as their group
of automorphism.

\bigskip Last, as a matter of notation, we will denote the norm on a Banach
space $E$ by $\left\Vert .\right\Vert _{E}$ where there is no ambiguity. The
norm $\left\Vert .\right\Vert _{\mathbb{C}^{n}}$ is the canonical hermitian norm.

\section{Spectra of Quotients of $F_{n}^{\infty}$}

\bigskip This section addresses the duality between the process of associating
an ideal of $F_{n}^{\infty}$ to a subset of the \emph{open unit ball}
$\mathbb{B}_{n}$ of $\mathbb{C}^{n}$\ and the computation of the spectrum of a
quotient of $F_{n}^{\infty}$. To this end, we first observe that since
$F_{n}^{\infty}$ is a dual space, given any weak* closed ideal $\mathcal{J}$
in $F_{n}^{\infty}$, the space $F_{n}^{\infty}/\mathcal{J}$ is dual to the
polar $\mathcal{J}^{\circ}$ of $\mathcal{J}$ in the predual of $F_{n}^{\infty
}$, and thus can be endowed with the weak*\ topology. We can thus define:

\begin{definition}
Let $\mathcal{J}$ be a weak*\ closed two-sided ideal in $F_{n}^{\infty}$. A
weak* continuous unital algebra homomorphism from $F_{n}^{\infty}/\mathcal{J}$
to $\mathbb{C}$ is called a weak* scalar representation of $F_{n}^{\infty
}/\mathcal{J}$, or when no confusion may arise, a representation of
$F_{n}^{\infty}/\mathcal{J}$.
\end{definition}

\begin{definition}
Let $\mathcal{J}$ be a weak*\ closed two-sided ideal in $F_{n}^{\infty}$. We
define the spectrum $\Sigma\left(  F_{n}^{\infty}/\mathcal{J}\right)  $ of
$F_{n}^{\infty}/\mathcal{J}$ as the set of all weak* scalar representations of
$F_{n}^{\infty}/\mathcal{J}$, endowed with the weak*\ topology of the dual of
$F_{n}^{\infty}/\mathcal{J}$.
\end{definition}

\bigskip Let us note that this definition involves two distinct topologies. As
a set, $\Sigma\left(  F_{n}^{\infty}/\mathcal{J}\right)  $ is defined as the
collection of scalar valued unital algebra homorphisms of $F_{n}^{\infty}$
which are continuous for the weak* topology of $F_{n}^{\infty}/\mathcal{J}$
seen as the dual space of $\mathcal{J}^{\circ}$. On the other hand, the
topology on $\Sigma\left(  F_{n}^{\infty}/\mathcal{J}\right)  $ is the
restriction of the weak* topology of the dual $\left(  F_{n}^{\infty
}/\mathcal{J}\right)  ^{\ast}$ of $F_{n}^{\infty}/\mathcal{J}$.

\bigskip Now, by \cite{Popescu91}, for any $\lambda=\left(  \lambda_{1}%
,\ldots,\lambda_{n}\right)  \in\mathbb{B}_{n}$ we can define a unique
representation $\pi_{\lambda}$ of $F_{n}^{\infty}$ such that $\pi_{\lambda
}\left(  S_{j}\right)  =\lambda_{j}$ with $j=1,\ldots,n$. Moreover, by
\cite{Popescu91} the map $\lambda\in\mathbb{B}_{n}\mapsto\pi_{\lambda}$ is an
homeomorphism. Given any $\varphi\in F_{n}^{\infty}$, we denote $\pi_{\lambda
}\left(  \varphi\right)  $ by $\varphi\left(  \lambda\right)  $, thus using
duality to see elements of $F_{n}^{\infty}$ as functions on $\mathbb{B}_{n}$.
With these notations, we have the following standard result:

\begin{proposition}
\label{trivial2}Let $\mathcal{J}$ be a weak*\ closed two-sided ideal in
$F_{n}^{\infty}$. Then:%
\[
\Sigma\left(  F_{n}^{\infty}/\mathcal{J}\right)  =\left\{  z\in\mathbb{B}%
_{n}:\forall\varphi\in\mathcal{J}\ \ \ \varphi\left(  z\right)  =0\right\}
\text{.}%
\]

\end{proposition}

\begin{proof}
The canonical projection $q:F_{n}^{\infty}\rightarrow F_{n}^{\infty
}/\mathcal{J}$ is weak* continuous, and thus if $\pi$ is a weak* continuous
algebra homomorphism from $F_{n}^{\infty}/\mathcal{J}$ to $\mathbb{C}$ then
$\pi\circ q$ is in the spectrum of $F_{n}^{\infty}$. Conversely, if $\pi$ is
in the spectrum of $F_{n}^{\infty}$ such that $\pi\left(  \mathcal{J}\right)
=\left\{  0\right\}  $ then $\pi$ defines an element of the spectrum of
$F_{n}^{\infty}/\mathcal{J}$. We thus have proven that the spectrum of
$F_{n}^{\infty}/\mathcal{J}$ is given by $\left\{  z\in\mathbb{B}_{n}%
:\forall\varphi\in F_{n}^{\infty}\ \ \ \varphi\left(  z\right)  =0\right\}  $.
\end{proof}

\bigskip We observe that, by Proposition (\ref{trivial2}), if $\varphi\in
F_{n}^{\infty}/\mathcal{J}$, where $\mathcal{J}$ is a weak* closed two-sided
ideal in $F_{n}^{\infty}$, and $\lambda\in\Sigma\left(  F_{n}^{\infty
}/\mathcal{J}\right)  $ then $\varphi\left(  \lambda\right)  $ is well-defined
as the common value of $\psi\left(  \lambda\right)  $ for $\psi\in
F_{n}^{\infty}$ such that $\psi+\mathcal{J}=\varphi$ (where $\psi+\mathcal{J}$
is the class of $\psi$ in $F_{n}^{\infty}/\mathcal{J}$). Now, we can use this
function representation of $F_{n}^{\infty}/\mathcal{J}$ on its spectrum to
associate a natural holomorphic map to any automorphism $\Phi$ of
$F_{n}^{\infty}/\mathcal{J}$. Indeed, if $\pi$ is a weak* scalar
representation of $F_{n}^{\infty}/\mathcal{J}$ then so is $\pi\circ\Phi$ for
all weak* continuous automorphism of $F_{n}^{\infty}/\mathcal{J}$. Hence, we
can define the following:

\begin{definition}
Let $\mathcal{J}$ be a weak*\ closed two-sided ideal in $F_{n}^{\infty}$. Let
$\Phi$ be a weak* continuous automorphism of $F_{n}^{\infty}/\mathcal{J}$. The
\emph{dual map} $\widehat{\Phi}$ of $\Phi$ is defined for all $\lambda
\in\Sigma\left(  F_{n}^{\infty}/\mathcal{J}\right)  $ by:%
\[
\widehat{\Phi}\left(  \lambda\right)  =\left(  \Phi^{-1}\left(  S_{1}\right)
\left(  \lambda\right)  ,\ldots,\Phi^{-1}\left(  S_{1}\right)  \left(
\lambda\right)  \right)  \text{.}%
\]

\end{definition}

\begin{definition}
The group of weak* completely isometric continuous automorphisms of
$F_{n}^{\infty}/\mathcal{J}$ is denoted by $\operatorname*{Aut}\left(
F_{n}^{\infty}/\mathcal{J}\right)  $.
\end{definition}

\bigskip It is an immediate observation that the map $\Phi\in
\operatorname*{Aut}\left(  F_{n}^{\infty}/\mathcal{J}\right)  \mapsto
\widehat{\Phi}$ is a group homomorphism to the group of homeomorphisms of
$\Sigma\left(  F_{n}^{\infty}/\mathcal{J}\right)  $ (in particular, mapping
the identity to the identity) such that:%
\[
\forall\varphi\in F_{n}^{\infty}/\mathcal{J}\ \ \ \forall\lambda\in
\Sigma\left(  F_{n}^{\infty}/\mathcal{J}\right)  \ \ \ \ \Phi(\varphi)\left(
\lambda\right)  =\varphi\left(  \widehat{\Phi}^{-1}\left(  \lambda\right)
\right)  \text{.}%
\]

\bigskip Although, for general completely contractive homomorphism $\Phi$, one
would be inclined to define the dual of $\Phi$ as $\Phi^{\prime}:\lambda
\in\Sigma\left(  F_{n}^{\infty}/\mathcal{J}\right)  \mapsto\left(  \Phi\left(
S_{j}\right)  \left(  \lambda\right)  \right)  _{j=1,\ldots,n}$, the resulting
map restricted to automorphisms would be valued in the opposite of the group
of homeomorphisms of $\Sigma\left(  F_{n}^{\infty}/\mathcal{J}\right)  $,
which will be inconvenient. Hence, we adopt our modified definition in this
paper, to obtain a group morphism. This is the same definition as in
\cite{Pitts98} when $\mathcal{J}=\left\{  0\right\}  $.

\bigskip Although a priori only a topological space homeomorphic to
$\mathbb{B}_{n}$, the spectrum $\Sigma\left(  F_{n}^{\infty}\right)  $ is in
fact endowed with a complex structure via its relation with $F_{n}^{\infty}$
\cite{Popescu91}: the maps $\varphi$ defined by elements of $F_{n}^{\infty} $
on $\mathbb{B}_{n}$ are holomorphic, and so are the dual maps of automorphisms
\cite[Theorem 4.11]{Pitts98}. These results can be extended to more general
complex domains \cite{Popescu08}. This allows us to use techniques from the
theory of analytic functions in several complex variables, as in
\cite{Arias09}.

\bigskip Our main focus in this section are the following two related notions:

\begin{definition}
Let $\Delta\subseteq\mathbb{B}_{n}$. The Nevanlinna ideal for $\Delta$ in
$F_{n}^{\infty}$\ is:%
\[
\mathcal{N}_{\Delta}=\left\{  \varphi\in F_{n}^{\infty}:\forall z\in
\Delta\ \ \ \varphi\left(  z\right)  =0\right\}
\]
and the quotient of $F_{n}^{\infty}$ localized at $\Delta$ is $F_{n,\Delta
}^{\infty}=F_{n}^{\infty}/\mathcal{N}_{\Delta}$.
\end{definition}

\begin{definition}
Let $\Delta\subseteq\mathbb{B}_{n}$. The spectrum of $F_{n,\Delta}^{\infty}$
is called the \emph{spectral closure} of $\Delta$ and is denoted by
$\overline{\Delta}^{\Sigma}$.
\end{definition}

\bigskip By Proposition (\ref{trivial2}), given $\Delta\subseteq\mathbb{B}%
_{n}$, we always have $\Delta\subseteq\overline{\Delta}^{\Sigma}$. In general,
we can have $\Delta\subsetneq\overline{\Delta}^{\Sigma}$, as shown for
instance when $\Delta=\left\{  0,\frac{1}{n}:n\in\mathbb{N}\text{,}%
n>0\right\}  \subset\mathbb{B}_{1}$, since there is no nonzero holomorphic
function which is null on $\Delta$, so $\overline{\Delta}^{\Sigma}%
=\mathbb{B}_{1}$. In other words, the main issue when relating $\Delta
\subseteq\mathbb{B}_{n}$ with $\overline{\Delta}^{\Sigma}$ is that the
Nevanlinna ideal $\mathcal{N}_{\Delta}$ may be null. It is well-known
\cite[9.1.4, 9.1.5]{Krantz99} that, for $n=1$, a set $\Delta$ is the zero set
for some holomorphic function if and only if it satisfies the Blaschke
condition, i.e. $\Delta=\left\{  \lambda_{j}:j\in\mathbb{N}\right\}  $ with
$\sum_{j=0}^{\infty}\left(  1-\left\vert \lambda_{j}\right\vert \right)
<\infty$. Under this condition, the Blachke product associated to $\Delta$ is
a holomorphic function which is zero exactly on $\Delta$. Unfortunately, such
a result does not hold in higher dimension \cite[Ch. 9]{Krantz}. However, we
shall now prove that the Blaschke condition is still sufficient to ensure that
$\Delta=\overline{\Delta}^{\Sigma}$ in $\mathbb{B}_{n}$.

To this end, we shall use the geometry of $\mathbb{B}_{n}$ by providing a
formula connecting the Poincare pseudohyperbolic metric on $\mathbb{B}_{n}$
with the unit ball of $F_{n}^{\infty}$. As a tool for our proof, we will use
the following lemma, which is a special case of \cite[Theorem 8.1.4]{Rudin81}
and which will find a role in the next section as well. We include a proof of
this lemma for the reader's convenience.

\begin{lemma}
\label{Schwarz}Let $\varphi$ be a holomorphic function from $\mathbb{B}_{n}$
to $\mathbb{B}_{k}$ for some nonzero natural $k$ and such that $\varphi(0)=0
$. Then for all $z\in\mathbb{B}_{n}$ we have $\left\Vert \varphi\left(
z\right)  \right\Vert _{\mathbb{C}^{k}}\leq\left\Vert z\right\Vert
_{\mathbb{C}^{n}}$.
\end{lemma}

\begin{proof}
Let $z\in\mathbb{B}_{n}$. Since the result is trivial for $z=0$, we shall
assume $z\not =0$. Let $\theta:\mathbb{C}^{k}\rightarrow\mathbb{C}$ be a
linear functional of norm $1$ (for the dual norm to the canonical Hermitian
norm on $\mathbb{C}^{k}$) such that $\left\vert \theta\circ\varphi\left(
z\right)  \right\vert =\left\Vert \varphi\left(  z\right)  \right\Vert
_{\mathbb{C}^{k}}$. We define the map $\varphi_{z}^{\theta}:\mathbb{B}%
_{1}\longrightarrow\mathbb{B}_{1}$ by $\varphi_{z}^{\theta}(t)=\theta
\circ\varphi\left(  t\frac{z}{\left\Vert z\right\Vert _{\mathbb{C}^{n}}%
}\right)  $ for $t\in\mathbb{B}_{1}$. By construction, $\varphi_{z}^{\theta} $
is holomorphic from the unit disk into itself and $\varphi_{z}^{\theta}\left(
0\right)  =0$. Hence, by the Schwarz lemma, we have for all $t\in
\mathbb{B}_{1}$ that $\left\vert \varphi_{z}^{\theta}(t)\right\vert
\leq\left\vert t\right\vert $. In particular:%
\[
\left\Vert z\right\Vert _{\mathbb{C}^{n}}\geq\left\vert \varphi_{z}^{\theta
}\left(  \left\Vert z\right\Vert \right)  \right\vert =\left\vert \theta
\circ\varphi\left(  z\right)  \right\vert =\left\Vert \varphi\left(  z\right)
\right\Vert _{\mathbb{C}^{k}}\text{.}%
\]
Hence our lemma is proven.
\end{proof}

\bigskip Poincare's pseudohyperbolic metric on the open unit ball
$\mathbb{B}_{n}$ of $\mathbb{C}^{n}$ between two points $z$ and $w$ can be
defined as the Euclidean distance $\rho$\ between $0$ and the image of $w$ by
any biholomorphic function of the ball which maps $z$ to $0$. As customary in
complex analysis, we will refer to biholomorphic maps of $\mathbb{B}_{n}$ onto
itself as \emph{automorphisms of }$\mathbb{B}_{n}$. We show that it is also
possible to compute this distance by using $F_{n}^{\infty}$.

\begin{proposition}
\label{metric}Let $\rho$ be the Poincare pseudohyperbolic metric on
$\mathbb{B}_{n}$. For any $z,w\in\mathbb{B}_{n}$ we have:%
\[
\rho\left(  z,w\right)  =\max\left\{  |\varphi(z)|:\varphi\in F_{n}^{\infty
}\text{ with }\left\Vert \varphi\right\Vert _{F_{n}^{\infty}}\leq1\text{ and
}\varphi\left(  w\right)  =0\right\}  \text{.}%
\]

\end{proposition}

\begin{proof}
We define for all $z,w\in\mathbb{B}_{n}$ the quantity
\[
\eta\left(  z,w\right)  =\sup\left\{  \left\vert \varphi(z)\right\vert
:\varphi\in F_{n}^{\infty}\text{ with }\left\Vert \varphi\right\Vert
_{F_{n}^{\infty}}\leq1\text{ and }\varphi\left(  w\right)  =0\right\}
\text{.}%
\]
We wish to show that $\eta=\rho$ and that the supremum defining $\eta$ is, in
fact, reached.

First, we prove that $\eta$ is invariant under the action of
$\operatorname*{SU}(n,1)$ on $\mathbb{B}_{n}$. Let $z,w\in\mathbb{B}_{n}$ and
let $\widehat{\Phi}$ be an automorphism of $\mathbb{B}_{n}$. There exists by
\cite{Pitts98}\ a unique automorphism $\Phi$ of $F_{n}^{\infty}$ such that,
for all $\omega\in\mathbb{B}_{n}$ and $\varphi\in F_{n}^{\infty}$, we have:
\[
\Phi^{-1}\left(  \varphi\right)  \left(  \omega\right)  =\varphi\left(
\widehat{\Phi}\left(  \omega\right)  \right)  \text{.}%
\]
(Of course, we could denote this automorphism of $F_{n}^{\infty}$ by $\Phi$
rather than $\Phi^{-1}$ but we prefer to keep the notations for dual map
consistent in this paper).

As an automorphism of $F_{n}^{\infty}$ is an isometry and thus maps the unit
ball of $F_{n}^{\infty}$ onto itself. Thus:%
\begin{align*}
\eta\left(  \widehat{\Phi}\left(  z\right)  ,\widehat{\Phi}\left(  w\right)
\right)   &  =\sup\left\{  \left\vert \varphi\left(  \widehat{\Phi}\left(
z\right)  \right)  \right\vert \left\vert
\begin{array}
[c]{l}%
\varphi\in F_{n}^{\infty}\text{,}\\
\left\Vert \varphi\right\Vert _{F_{n}^{\infty}}\leq1\text{,}\\
\varphi\left(  \widehat{\Phi}\left(  w\right)  \right)  =0
\end{array}
\right.  \right\} \\
&  =\sup\left\{  \left\vert \Phi^{-1}\left(  \varphi\right)  \left(  z\right)
\right\vert \left\vert
\begin{array}
[c]{l}%
\varphi\in F_{n}^{\infty}\text{,}\\
\left\Vert \Phi^{-1}\left(  \varphi\right)  \right\Vert _{F_{n}^{\infty}}%
\leq1\text{,}\\
\Phi^{-1}\left(  \varphi\right)  \left(  w\right)  =0
\end{array}
\right.  \right\} \\
&  =\sup\left\{  \left\vert \psi\left(  z\right)  \right\vert \left\vert
\begin{array}
[c]{l}%
\psi\in F_{n}^{\infty}\text{,}\\
\left\Vert \psi\right\Vert _{F_{n}^{\infty}}\leq1\text{,}\\
\psi\left(  w\right)  =0
\end{array}
\right.  \right\}  =\eta\left(  z,w\right)  \text{.}%
\end{align*}
In particular, $\eta\left(  z,w\right)  =\eta\left(  w,z\right)  $ as there
exists an automorphism of $\mathbb{B}_{n}$ which maps $z$ to $w$ and
vice-versa. Thus, it is enough to prove that, for any $z\in\mathbb{B}_{n}$, we
have $\rho\left(  0,z\right)  =\left\Vert z\right\Vert _{\mathbb{C}^{n}}$.
This would suffice to show that $\eta$ is the Poincare pseudohyperbolic metric
$\rho$ on $\mathbb{B}_{n}$.

Let us fix $z\in\mathbb{B}_{n}$. Since $\eta\left(  0,0\right)  =0$, we may as
well assume $z\not =0$. Let $\varphi\in F_{n}^{\infty}$ such that $\left\Vert
\varphi\right\Vert _{F_{n}^{\infty}}\leq1$ and $\varphi\left(  0\right)  =0$.
By Lemma (\ref{Schwarz}), we have $\left\vert \varphi\left(  z\right)
\right\vert \leq\left\Vert z\right\Vert _{\mathbb{C}^{n}}$ so $\eta\left(
0,z\right)  $, which is the supremum of $\left\vert \varphi\left(  z\right)
\right\vert $ for $\varphi\in F_{n}^{\infty}$ with $\left\Vert \varphi
\right\Vert _{F_{n}^{\infty}}\leq1$ and $\varphi(0)=0$, is bounded above by
$\left\Vert z\right\Vert _{\mathbb{C}^{n}}$. On the other hand, observe that
for any $a_{1},\ldots,a_{n}\in\mathbb{C}^{n}$, if $\varphi=\sum_{j=1}^{n}%
a_{j}S_{j}$ then the norm of $\varphi$ is $\left\Vert \varphi^{\ast}%
\varphi\right\Vert _{\mathcal{B}\left(  F_{n}^{2}\right)  }^{\frac{1}{2}}$
which equals $\sqrt[2]{\sum_{j=1}^{n}\left\vert a_{j}\right\vert ^{2}}$ since
$S_{j}^{\ast}S_{k}=\delta_{j}^{k}1$. If we wri$_{{}}$te $z=\left(
z_{1},\ldots,z_{n}\right)  $ then, choosing $\varphi_{z}=\sum_{j=1}^{n}%
\frac{z_{i}}{\left\Vert z\right\Vert _{\mathbb{C}^{n}}}S_{i}$ we see that
$\varphi_{z}\in F_{n}^{\infty}$ with $\varphi_{z}(0)=0$ and $\left\vert
\varphi_{z}\left(  z\right)  \right\vert =\left\Vert z\right\Vert
_{\mathbb{C}^{n}}$. So $\left\Vert z\right\Vert _{\mathbb{C}^{n}}\leq
\eta\left(  0,z\right)  $ as desired. We conclude that $\eta(0,z)=\left\Vert
z\right\Vert _{\mathbb{C}^{n}}$ and this supremum is reached at $\varphi_{z}$.
\end{proof}

\bigskip We now can prove that the Blaschke condition is sufficient for a
subset $\Delta$ of $\mathbb{B}_{n}$ to equal its spectral closure. We start
with the following lemma which uses an important estimate from the theory of
functions on $\mathbb{B}_{n}$.

\begin{lemma}
\label{blaschke_lemma}Let $\left\{  \lambda_{j}:j\in\mathbb{N}\right\}
\subseteq\mathbb{B}_{n}$. Let $\varphi$ be an automorphism of $\mathbb{B}_{n}$
such that $\varphi\left(  0\right)  \not =0$. Then:%
\[
\sum_{j=0}^{\infty}\left(  1-\left\Vert \lambda_{j}\right\Vert _{\mathbb{C}%
^{n}}\right)  <\infty\iff\sum_{j=0}^{\infty}\left(  1-\left\Vert
\varphi\left(  \lambda_{j}\right)  \right\Vert _{\mathbb{C}^{n}}\right)
<\infty\text{.}%
\]

\end{lemma}

\begin{proof}
We denote $\left\Vert .\right\Vert _{\mathbb{C}^{n}}$ by $\left\Vert
.\right\Vert $ in this proof. Let $a=\varphi^{-1}(0)$ and note that $a\not =0
$ by assumption. Using \cite[Theorem 2.2.2 p. 26]{Rudin81}, we have for all
$z\in\mathbb{B}_{n}$:%
\[
1-\left\Vert \varphi\left(  z\right)  \right\Vert =\frac{\left(  1-\left\Vert
a\right\Vert ^{2}\right)  \left(  1+\left\Vert z\right\Vert \right)  }{\left(
1+\left\Vert \varphi\left(  z\right)  \right\Vert \right)  \left\vert
1-\left\langle a,z\right\rangle \right\vert }\left(  1-\left\Vert z\right\Vert
\right)
\]
where $\left\langle .,.\right\rangle $ is the canonical inner product in
$\mathbb{C}^{n}$. Since, for all $z\in\mathbb{B}_{n}$, we have:%
\[
0<\frac{1-\left\Vert a\right\Vert ^{2}}{2}\leq\frac{\left(  1-\left\Vert
a\right\Vert ^{2}\right)  \left(  1+\left\Vert z\right\Vert \right)  }{\left(
1+\left\Vert \varphi\left(  z\right)  \right\Vert \right)  \left\vert
1-\left\langle a,z\right\rangle \right\vert }\leq\frac{2}{1-\left\Vert
a\right\Vert }\text{.}%
\]
The result follows.
\end{proof}

\begin{theorem}
\label{Blaschke}Let $\Delta=\left\{  \lambda_{j}:j\in\mathbb{N}\right\}
\subseteq\mathbb{B}_{n}$. If $\sum_{j=0}^{\infty}\left(  1-\left\Vert
\lambda_{j}\right\Vert _{\mathbb{C}^{n}}\right)  <\infty$ then the weak*
spectrum $\Sigma\left(  F_{n,\Delta}^{\infty}\right)  $ of $F_{n,\Delta
}^{\infty}=F_{n}^{\infty}/\mathcal{N}_{\Delta}$ is $\Delta$.
\end{theorem}

\begin{proof}
Assume given $\Delta=\left\{  \lambda_{j}:j\in\mathbb{N}\right\}
\subseteq\mathbb{B}_{n}$ with $\sum_{j=0}^{\infty}\left(  1-\left\Vert
\lambda_{j}\right\Vert _{\mathbb{C}^{n}}\right)  <\infty$. Let $z\in
\mathbb{B}_{n}\setminus\Delta$. We wish to show that there exists $\varphi
\in\mathcal{N}_{\Delta}$ such that $\varphi\left(  z\right)  \not =0$. Let
$\widehat{\Phi}$ be an automorphism of $\mathbb{B}_{n}$ which maps $z$ to $0$.
By Lemma (\ref{metric}), for each $j\in\mathbb{N}$ there exists $\varphi
_{j}\in F_{n}^{\infty}$ such that $\left\Vert \varphi_{j}\right\Vert
_{\mathbb{C}^{n}}\leq1$, $\varphi_{j}(\widehat{\Phi}(\lambda_{j}))=0$ and
$\left\Vert \widehat{\Phi}(\lambda_{j})\right\Vert _{\mathbb{C}^{n}}%
=\varphi_{j}\left(  0\right)  $. Fix $j\in\mathbb{N}$. We define the element
$\psi_{j}=%
{\displaystyle\prod\limits_{k=0}^{j}}
\varphi_{k}\in F_{n}^{\infty}$. By construction, we have $\left\Vert \psi
_{j}\right\Vert _{F_{n}^{\infty}}\leq1$, as well as $\psi_{j}\left(  0\right)
=%
{\displaystyle\prod\limits_{k=0}^{j}}
\left\Vert \widehat{\Phi}(\lambda_{j})\right\Vert _{\mathbb{C}^{n}}$ and
$\psi_{j}\left(  \widehat{\Phi}(\lambda_{k})\right)  =0$ for $k\in\left\{
0,\ldots,j\right\}  $. Now, the unit ball of $F_{n}^{\infty}$ is weak*
compact, so we can extract a subsequence of $\left(  \psi_{j}\right)
_{j\in\mathbb{N}}$ which converges in the weak* topology to some $\psi\in
F_{n}^{\infty}$. We recall that by definition, the notation $\psi(\mu)$ refers
to $\pi_{\mu}\left(  \psi\right)  $ where $\pi_{\mu}$ is the unique unital
weak* continuous algebra homomorphism from $F_{n}^{\infty}$ to $\mathbb{C}$
mapping the canonical generators $S_{1},\ldots,S_{n}$ of $F_{n}^{\infty}$ to
$\mu_{1},\ldots,\mu_{n}$ with $\mu=\left(  \mu_{1},\ldots,\mu_{n}\right)  $.
Hence, by continuity, we have in particular, since $\sum_{j=0}^{\infty}\left(
1-\left\Vert \lambda_{j}\right\Vert _{\mathbb{C}^{n}}\right)  <\infty$, using
Lemma (\ref{blaschke_lemma}), we have $\sum_{j=0}^{\infty}\left(  1-\left\Vert
\widehat{\Phi}(\lambda_{j})\right\Vert _{\mathbb{C}^{n}}\right)  <\infty$ and
thus:%
\[
\psi\left(  0\right)  =%
{\displaystyle\prod\limits_{j=0}^{\infty}}
\left\Vert \widehat{\Phi}(\lambda_{j})\right\Vert _{\mathbb{C}^{n}}>0\text{.}%
\]
Hence $\psi(0)\not =0$ while $\psi\left(  \widehat{\Phi}(\lambda_{j})\right)
=0$ by construction as well. Thus, if $\Phi$ is the automorphism of
$F_{n}^{\infty}$ whose dual map is $\widehat{\Phi}$ and if we set
$\varphi=\Phi^{-1}\left(  \psi\right)  $, we see that $\varphi\in
\mathcal{N}_{\Delta}$ while $\varphi(z)\not =0$ so $z\not \in \overline
{\Delta}^{\Sigma}$. Thus $\overline{\Delta}^{\Sigma}\subseteq\Delta$. Since
the reverse inclusion is Proposition (\ref{trivial2}), our theorem is established.
\end{proof}

\section{Automorphism groups of quotients of $F_{n}^{\infty}$}

\bigskip This section establishes that, under a natural condition of
convergence, it is possible to choose many discrete subgroups of
$\operatorname*{SU}(n,1)$ as full automorphism groups of some operator
algebras obtained as quotients of $F_{n}^{\infty}$. We shall call any
biholomorphic from $\mathbb{B}_{n}$ onto $\mathbb{B}_{n}$ an
\emph{automorphism} of $\mathbb{B}_{n}$. The group of automorphisms of
$\mathbb{B}_{n}$ is denoted by $\operatorname*{Aut}\left(  \mathbb{B}%
_{n}\right)  $.

\bigskip We start with an easy consequence of \cite{Pitts98}:

\begin{lemma}
\label{trivial3}Let $\Delta\subseteq\mathbb{B}_{n}$. Let $\phi$ be an
automorphism of $\mathbb{B}_{n}$ such that $\phi\left(  \overline{\Delta
}^{\Sigma}\right)  \subseteq\overline{\Delta}^{\Sigma}$. Then there exists an
automorphism $\Phi$ of $F_{n}^{\infty}$ such that $\widehat{\Phi}=\phi^{-1}$
and $\Phi\left(  \mathcal{N}_{\Delta}\right)  \subseteq\mathcal{N}_{\Delta}$,
so that $\Phi$ induces an automorphism of $F_{n}^{\infty}/\mathcal{N}_{\Delta
}$.
\end{lemma}

\begin{proof}
If $\phi$ is any automorphism of $\mathbb{B}_{n}$ such that $\phi\left(
\overline{\Delta}^{\Sigma}\right)  \subseteq\overline{\Delta}^{\Sigma}$, then
by \cite{Pitts98} there exists a unique automoprhism $\Phi$ of $F_{n}^{\infty
}$ such that $\phi^{-1}$ is the dual map of $\Phi$ on the spectrum
$\mathbb{B}_{n}$ of $F_{n}^{\infty}$. Now, let $x\in\mathcal{N}_{\Delta}$. By
construction, if $\lambda\in\overline{\Delta}^{\Sigma}$ then $\Phi\left(
x\right)  \left(  \lambda\right)  =x\left(  \phi\left(  \lambda\right)
\right)  =0$ since $\phi\left(  \lambda\right)  \in\overline{\Delta}^{\Sigma}%
$. Hence $\Phi\left(  \mathcal{N}_{\Delta}\right)  \subseteq\mathcal{N}%
_{\Delta}$ and thus $\Phi$ induces an automorphism of $F_{n,\Delta}^{\infty}$.
\end{proof}

\bigskip The converse implication, i.e. that any automorphism of
$F_{n}^{\infty}/\mathcal{N}_{\Delta}$ is given by an automorphism of
$\mathbb{B}_{n}$ which maps the spectrum to itself, is the subject of the rest
of this paper. To this end, we shall use the following special case of
\cite[Theorem 8.2.2]{Rudin81}:

\begin{lemma}
\label{Schwarz2}Let $\Phi:\mathbb{B}_{n}\rightarrow\mathbb{B}_{n}$ be a
holomorphic map with $\Phi(0)=0$. If the set $\left\{  \lambda\in
\mathbb{B}_{n}:\Phi(\lambda)=\lambda\right\}  $ of fixed points of $\Phi$
spans $\mathbb{C}^{n}$, then $\Phi$ is the identity.
\end{lemma}

\begin{proof}
Our assumption implies the existence of a basis of invariant vectors for the
Frechet derivative $\Phi^{\prime}(0)$ of $\Phi$ at $0$, so $\Phi^{\prime}(0)$
is the identity on $\mathbb{C}^{n}$. By \cite[Theorem 8.2.2]{Rudin81}, $\Phi$
has the same invariant points as $\Phi^{\prime}(0)$, hence our lemma.
\end{proof}

\bigskip We now establish the key result for this section. The class of $x\in
F_{n}^{\infty}$ in $F_{n,\Delta}^{\infty}$ is denoted by $x+\mathcal{N}%
_{\Delta}$. To make the proof of Theorem (\ref{main}) clearer, we organize it
as a succession of lemmas. The first step is to construct automorphisms of
$\mathbb{B}_{n}$ from automorphisms of $F_{n,\Delta}^{\infty}$ ($\Delta
\subseteq\mathbb{B}_{n}$).

\begin{lemma}
\bigskip\label{construction}Let $\Delta\subseteq\mathbb{B}_{n}$ such that
$\overline{\Delta}^{\Sigma}$ spans $\mathbb{C}^{n}$. Let $\Phi$ be a
completely contractive automorphism of $F_{n,\Delta}^{\infty}$. Then there
exists $\varphi_{1},\ldots,\varphi_{n}\in F_{n}^{\infty}$ such that:

\begin{itemize}
\item We have $\left\Vert \left[  \varphi_{1},\ldots,\varphi_{n}\right]
\right\Vert _{M_{1,n}\left(  F_{n}^{\infty}\right)  }\leq1$,

\item We have $\left[  \Phi^{-1}\left(  S_{j}+\mathcal{N}_{\Delta}\right)
\right]  _{j=1,\ldots,n}=\left[  \varphi_{j}\right]  _{j=1,\ldots,n}%
+M_{1,n}\left(  \mathcal{N}_{\Delta}\right)  $,

\item The dual map $\widehat{\Phi}$ of $\Phi$ on the spectrum $\overline
{\Delta}^{\Sigma}$ of $F_{n,\Delta}^{\infty}$ is given by $\lambda
\mapsto\left(  \varphi_{1}\left(  \lambda\right)  ,\ldots,\varphi_{n}%
(\lambda)\right)  $,

\item The map $\overleftrightarrow{\Phi}=\left(  \varphi_{1},\ldots
,\varphi_{n}\right)  $ is a biholomorphic map from $\mathbb{B}_{n}$ onto
$\mathbb{B}_{n}$.
\end{itemize}
\end{lemma}

\begin{proof}
Let $\Phi$ be an automorphism of $F_{n,\Delta}^{\infty}$. The element:%
\[
\left[  \Phi^{-1}(S_{1}+\mathcal{N}_{\Delta})\ \cdots\ \Phi^{-1}\left(
S_{n}+\mathcal{N}_{\Delta}\right)  \right]  \in M_{1,n}\left(  F_{n,\Delta
}^{\infty}\right)  =M_{1,n}\left(  F_{n}^{\infty}\right)  /M_{1,n}\left(
\mathcal{N}_{\Delta}\right)
\]
is a row contraction, since $\left[  S_{1}\ \cdots\ S_{n}\right]  $ is and
$\Phi$ is completely isometric, hence so is $\Phi^{-1}$. Now, by definition:%
\begin{align*}
1  &  \geq\left\Vert \left[  \Phi^{-1}(S_{1}+\mathcal{N}_{\Delta}%
)\ \cdots\ \Phi^{-1}\left(  S_{n}+\mathcal{N}_{\Delta}\right)  \right]
\right\Vert _{M_{1,n}\left(  F_{n,\Delta}^{\infty}\right)  }\\
&  =\inf_{\substack{\psi_{j}+\mathcal{N}_{\Delta}=\Phi^{-1}\left(
S_{j}\right)  \\j=1,\ldots,n}}\left\Vert \left[  \psi_{1}\ \cdots\ \psi
_{n}\right]  \right\Vert _{M_{1,n}\left(  F_{n}^{\infty}\right)  }\text{.}%
\end{align*}

Since the unit ball of $F_{n}^{\infty}$ is weak* compact, and hence so is the
unit ball of $M_{1,n}\left(  F_{n}^{\infty}\right)  $, and since
$\mathcal{N}_{\Delta}$ is weak*\ closed, we can find $\varphi_{1}%
,\ldots,\varphi_{n}\in F_{n}^{\infty}$ such that $\left\Vert \left[
\varphi_{1}\ \cdots\ \varphi_{n}\right]  \right\Vert _{M_{1,n}}=1$ and:%
\[
\left[  \Phi^{-1}(S_{1}+\mathcal{N}_{\Delta})\ \cdots\ \Phi^{-1}\left(
S_{n}+\mathcal{N}_{\Delta}\right)  \right]  =\left[  \varphi_{1}%
\ \cdots\ \varphi_{n}\right]  +M_{1,n}\left(  \mathcal{N}_{\Delta}\right)
\text{.}%
\]

Now, let $\widehat{\Phi}$ be the dual map of $\Phi$ on $\overline{\Delta
}^{\Sigma}$. Note that in particular $\widehat{\Phi}\left(  \overline{\Delta
}^{\Sigma}\right)  \subseteq\overline{\Delta}^{\Sigma}\subseteq\mathbb{B}_{n}%
$. Let $\lambda\in\overline{\Delta}^{\Sigma}$. Then:%
\begin{align*}
\widehat{\Phi}\left(  \lambda\right)   &  =\left(  \left[  \Phi^{-1}\left(
S_{1}+\mathcal{N}_{\Delta}\right)  \right]  \left(  \lambda\right)
,\ldots,\left[  \Phi^{-1}\left(  S_{n}+\mathcal{N}_{\Delta}\right)  \right]
\left(  \lambda\right)  \right) \\
&  =\left(  \left(  \varphi_{1}+\mathcal{N}_{\Delta}\right)  \left(
\lambda\right)  ,\ldots,\varphi_{n}+\mathcal{N}_{\Delta}\left(  \lambda
\right)  \right) \\
&  =\left(  \varphi_{1}\left(  \lambda\right)  ,\ldots,\varphi_{n}\left(
\lambda\right)  \right)  \text{ since for all }\theta\in\mathcal{N}_{\Delta
}\text{ we have }\theta\left(  \lambda\right)  =0\text{.}%
\end{align*}

Thus $\overleftrightarrow{\Phi}:=\left(  \varphi_{1},\ldots,\varphi
_{n}\right)  $, which is a holomorphic map from $\mathbb{B}_{n}$ to
$\overline{\mathbb{B}_{n}}$ (since $\left\Vert \left[  \varphi_{1}%
~\cdots~\varphi_{n}\right]  \right\Vert \leq1$), is an analytic extension of
$\widehat{\Phi}$ to $\mathbb{B}_{n}$. Now, let $z\in\mathbb{B}_{n}$ and
suppose that $\overleftrightarrow{\Phi}\left(  z\right)  $ lies on the
boundary of $\overline{\mathbb{B}_{n}}$. Then, up to conjugating
$\overleftrightarrow{\Phi}$ by a biholomorphic map, we may as well assume that
$\overleftrightarrow{\Phi}\left(  0\right)  $ lies on the boundary of
$\overline{\mathbb{B}_{n}}$. Again, up to conjugation by a unitary, we may as
well assume that $\overleftrightarrow{\Phi}\left(  0\right)  =\left(
1,0,\ldots,0\right)  $. Thus $\varphi_{1}(0)=1$ and $\varphi_{1}%
:\mathbb{B}_{1}\rightarrow\overline{\mathbb{B}_{1}}$ is holomorphic. We
conclude by the maximum modulus principle that $\varphi_{1}$ is the constant
function $1 $ on $\mathbb{B}_{1}$. Therefore $\varphi_{2}=\ldots=\varphi
_{n}=0$. Therefore, $\overleftrightarrow{\Phi}$ maps all of $\mathbb{B}_{n}$
to a constant value on the boundary of $\overline{\mathbb{B}_{n}}$, which
contradicts the fact that $\overleftrightarrow{\Phi}\left(  \overline{\Delta
}^{\Sigma}\right)  \subseteq\mathbb{B}_{n}$. Hence, $\overleftrightarrow{\Phi
}$ is a holomorphic map from $\mathbb{B}_{n}$ into $\mathbb{B}_{n}$.

With the same technique, we can construct a holomorphic map
$\overleftrightarrow{\Phi^{-1}}$ from $\mathbb{B}_{n}$ into itself whose
restriction to $\overline{\Delta}^{\Sigma}$ is the map $\widehat{\Phi^{-1}}$
dual to the inverse $\Phi^{-1}$ of $\Phi$ on $\overline{\Delta}^{\Sigma}$.

Now, by construction, $\overleftrightarrow{\Phi^{-1}}\circ\overleftrightarrow
{\Phi}\left(  \lambda\right)  =\lambda$ for all $\lambda\in$ $\overline
{\Delta}^{\Sigma}$ and $\overleftrightarrow{\Phi^{-1}}\circ\overleftrightarrow
{\Phi}$ is a holomorphic function from $\mathbb{B}_{n}$ to $\mathbb{B}_{n}$.
Since the span of $\overline{\Delta}^{\Sigma}$ is $\mathbb{C}^{n}$, we
conclude with Lemma (\ref{Schwarz2}) that $\overleftrightarrow{\Phi^{-1}}%
\circ\overleftrightarrow{\Phi}$ is the identity of $\mathbb{B}_{n}$. The same
exact reasoning shows that $\overleftrightarrow{\Phi}\circ\overleftrightarrow
{\Phi^{-1}}$ is also the identity of $\mathbb{B}_{n}$.

This concludes our lemma.
\end{proof}

Our next lemma establishes that the map $\Phi\in\operatorname*{Aut}\left(
F_{n,\Delta}^{\infty}\right)  \mapsto\widehat{\Phi}$ which maps an
automorphism of $F_{n,\Delta}^{\infty}$ to its dual map on the spectrum of
$F_{n,\Delta}^{\infty}$ is in fact a group monomorphism when $\Delta$ spans
$\mathbb{C}^{n}$. We briefly recall from \cite{Popescu91} the following
construction. We abbreviate the notation $F^{2}\left(  \mathbb{C}_{n}\right)
$ into $F_{n}^{2}$. Given a vector $\xi\in F_{n}^{2}$, we define $\left\Vert
\xi\right\Vert $ as the supremum of $\left\Vert \xi\otimes\eta\right\Vert
_{F_{n}^{2}}$ over finite linear combination of elementary tensors $\eta$ in
$F_{n}^{2}$ with $\left\Vert \eta\right\Vert _{F_{n}^{2}}\leq1$. If $\xi\in
F_{n}^{2}$ and $\left\Vert \xi\right\Vert <\infty$ then the operator $\eta\in
F_{n}^{2}\mapsto\xi\otimes\eta$, still denoted by $\xi$, is a well-defined
linear operator of norm $\left\Vert \xi\right\Vert $. From \cite{Popescu91},
we see that $F_{n}^{\infty}=\left\{  \xi\in F_{n}^{2}:\left\Vert
\xi\right\Vert <\infty\right\}  $ and $\left\Vert \cdot\right\Vert =\left\Vert
\cdot\right\Vert _{F_{n}^{\infty}}$. We also note that $\left\Vert
\cdot\right\Vert _{F_{n}^{\infty}}\geq\left\Vert \cdot\right\Vert _{F_{n}^{2}%
}$. With this identification in mind, we show:

\begin{lemma}
\label{uniqueness}Let $\Delta\subseteq\mathbb{B}_{n}$ such that $\Delta$ spans
$\mathbb{C}^{n}$. Let $\Phi$ be a completely contractive automorphism of
$F_{n,\Delta}^{\infty}$. If the dual map $\widehat{\Phi}$ of $\Phi$ on the
spectrum $\overline{\Delta}^{\Sigma}$ of $F_{n,\Delta}^{\infty}$ is the
identity, then $\Phi$ is the identity.

Consequently, $\Phi\in\operatorname*{Aut}\left(  F_{n,\Delta}^{\infty}\right)
\mapsto\widehat{\Phi}$ which maps an automorphism of $F_{n,\Delta}^{\infty}$
to its dual map is a group monomorphism.
\end{lemma}

\begin{proof}
By Lemma\ (\ref{construction}) there exists $\varphi_{1},\ldots,\varphi_{n}\in
F_{n}^{\infty}$ such that:%
\[
\left[  \Phi^{-1}\left(  S_{j}+\mathcal{N}_{\Delta}\right)  \right]
_{j=1\ldots n}=\left[  \varphi_{j}\right]  _{j=1\ldots n}+M_{1,n}\left(
\mathcal{N}_{\Delta}\right)
\]
and $\left\Vert \left[  \varphi_{j}\right]  _{j=1\ldots n}\right\Vert
_{M_{1,n}}\leq1$. Now, $\overleftrightarrow{\Phi}=\left(  \varphi_{1}%
,\ldots,\varphi_{n}\right)  $ fixes $\overline{\Delta}^{\Sigma}$ by
assumption, so $\overleftrightarrow{\Phi}$ is the identity by Lemma
(\ref{Schwarz2}). So $\varphi_{j}=S_{j}+C_{j}$ for some $C_{j}\in
F_{n}^{\infty} $ such that $C_{j}\left(  \lambda\right)  =0$ for all
$\lambda\in\mathbb{B}_{n}$, where $j=1,\ldots,n$. Hence $C_{j}$ lies in the
ideal of $F_{n}^{\infty}$ generated by $S_{j}S_{k}-S_{k}S_{j}$ ($k,j=1,\ldots
,n$).

Let $j\in\left\{  1,\ldots,n\right\}  $. Let $e_{j}$ be the $j^{\text{th}}$
canonical basis vector in $\mathbb{C}^{n}$ and let $\xi_{j}\in F_{n}^{2}$ such
that $C_{j}\left(  \eta\right)  =\xi_{j}\otimes\eta$ for all $\eta\in
F_{2}^{n}$. Then $\xi_{j}$ is orthogonal to $e_{j}$ in $F_{2}^{n}$ and
therefore:%
\[
1\geq\left\Vert \varphi_{j}\right\Vert _{F_{n}^{\infty}}^{2}\geq\left\Vert
e_{j}+\xi_{j}\right\Vert _{F_{n}^{2}}^{2}=\sqrt[2]{\left\Vert e_{j}\right\Vert
_{F_{n}^{2}}+\left\Vert \xi_{j}\right\Vert _{F_{n}^{2}}^{2}}\text{.}%
\]
Since $\left\Vert e_{j}\right\Vert =1$ we conclude $\xi_{j}=0$ and thus
$C_{j}=0$. Hence $\Phi^{-1}(S_{j}+\mathcal{N}_{\Delta})=S_{j}+\mathcal{N}%
_{\Delta}$ so $\Phi$ is the identity on $F_{n.\Delta}^{\infty}$ as desired.

In particular, the kernel of $\Phi\in\operatorname*{Aut}\left(  F_{n,\Delta
}^{\infty}\right)  \mapsto\widehat{\Phi}$ is reduced to the identity, so this
group homomorphism is injective.
\end{proof}

\begin{theorem}
\label{Lift}Let $\Delta\subseteq\mathbb{B}_{n}$ and let $\Gamma$ be the
subgroup of $\operatorname*{SU}(n,1)$ such that $\gamma\in\Gamma$ if and only
if $\gamma\left(  \overline{\Delta}^{\Sigma}\right)  =\overline{\Delta
}^{\Sigma}$. If the span of $\overline{\Delta}^{\Sigma}$ is $\mathbb{C}^{n} $
then the automorphism group of $F_{n,\Delta}^{\infty}$ is $\Gamma$.
\end{theorem}

\begin{proof}
Let $\Phi$ be an automorphism of $F_{n,\Delta}^{\infty}$. By Lemma
(\ref{construction}), we can extend the dual map $\widehat{\Phi}$ on the
spectrum $\overline{\Delta}^{\Sigma}$ of $F_{n,\Delta}^{\infty}$ into an
automorphism $\overleftrightarrow{\Phi}$ of $\mathbb{B}_{n}$ which maps
$\overline{\Delta}^{\Sigma}$ onto itself. By Lemma (\ref{uniqueness}), the map
$\Phi\in\operatorname*{Aut}\left(  F_{n,\Delta}^{\infty}\right)
\mapsto\widehat{\Phi}$ is a group monomorphism into the group of homeomorphism
of $\overline{\Delta}^{\Sigma}$. Now, by Lemma (\ref{Schwarz2}), the map
$\widehat{\Phi}\mapsto\overleftrightarrow{\Phi}$ constructed in
Lemma\ (\ref{construction}) is also a group monomorphism. Indeed, if any two
automorphisms of $\mathbb{B}_{n}$ agree on $\overline{\Delta}^{\Sigma}$ then
they must be equal. Since $\overleftrightarrow{\Phi_{1}}\circ
\overleftrightarrow{\Phi_{2}}$ certainly restricts on $\overline{\Delta
}^{\Sigma}$ to $\widehat{\Phi_{1}}\circ\widehat{\Phi_{2}}=\widehat{\Phi
_{1}\circ\Phi_{2}}$, we conclude by uniqueness that $\overleftrightarrow
{\Phi_{1}}\circ\overleftrightarrow{\Phi_{2}}=\overleftrightarrow{\Phi_{1}%
\circ\Phi_{2}}$ for any $\Phi_{1},\Phi_{2}\in\operatorname*{Aut}\left(
F_{n,\Delta}^{\infty}\right)  $. Hence, the map $\Phi\in\operatorname*{Aut}%
\left(  F_{n,\Delta}^{\infty}\right)  \mapsto\overleftrightarrow{\Phi}%
\in\operatorname*{Aut}\left(  \mathbb{B}_{n}\right)  $ is a group
monomorphism. Moreover, by Lemma (\ref{construction}), its range is included
in the set of automorphisms of $\mathbb{B}_{n}$ which maps $\overline{\Delta
}^{\Sigma}$ to itself. The reverse inclusion is established by Lemma
(\ref{trivial3}).

Thus, the map $\Phi\in\operatorname*{Aut}\left(  F_{n,\Delta}^{\infty}\right)
\longmapsto\overleftrightarrow{\Phi}\in\operatorname*{Aut}\left(
\mathbb{B}_{n}\right)  $ is a well-defined group isomorphism. This concludes
the computation of the automorphism group of $F_{n,\Delta}^{\infty}$.
\end{proof}

\begin{remark}
\bigskip As a by-product of the proof above, we see that any automorphism of
$F_{n,\Delta}^{\infty}$ lifts uniquely to an automorphism of $F_{n}^{\infty}$.
Even more: the group of automorphisms of $F_{n,\Delta}^{\infty}$ is the
quotient of the group of automorphisms of $F_{n}^{\infty}$ by the stabilizer
subgroup of $\mathcal{N}_{\Delta}$.
\end{remark}

\bigskip We conclude this section with the main result of this paper. Given a
discrete subgroup $\Gamma$ of $\operatorname*{SU}(n,1)$, when is the
automorphism group of a quotient of $F_{n}^{\infty}$ isomorphic to $\Gamma$?
In general, it is not true that the automorphism group of $F_{n,\Gamma
(0)}^{\infty}$ is $\Gamma$, as shown in the following example.

\begin{example}
\label{counter1}Let $\Gamma=\left\{  \gamma^{n}:n\in\mathbb{Z}\right\}  $
where $\gamma:z\in\mathbb{B}_{1}\mapsto\frac{z-\frac{1}{2}}{1-\frac{1}{2}z}$.
Then the rotation of center $0$ and angle $\pi$ is an (elliptic)\ automorphism
of $\mathbb{B}_{n}$ which maps the orbit $\Gamma(0)$ of $0$ for $\Gamma$ to
itself. Thus by Proposition (\ref{Lift}), it is an automorphism of
$F_{n,\Delta}^{\infty}$. Yet it is not an element of $\Gamma$, which only
consists of hyperbolic automorphisms.
\end{example}

\bigskip In sight of Example (\ref{counter1}), it is only natural to define:

\begin{definition}
Let $\Gamma$ be a subgroup of $\operatorname*{SU}\left(  n,1\right)  $. Let
$\underline{\Gamma}$ be the stabilizer subgroup of the orbit $\Gamma(0)$ of
$0$ by $\Gamma$ in $\operatorname*{SU}(n,1)$, i.e.:%
\[
\underline{\Gamma}=\left\{  \gamma\in\operatorname*{SU}\left(  n,1\right)
:\forall z\in\Gamma(0)\ \ \ \gamma\left(  z\right)  \in\Gamma\left(  0\right)
\right\}  \text{.}%
\]

\end{definition}

Example (\ref{counter1}) shows that $\Gamma$ may be a strict subgroup of
$\underline{\Gamma}$.

\bigskip

We can now put our results together to show that:

\begin{theorem}
\label{main}Let $\Gamma$ be a discrete subgroup of $\operatorname*{SU}(n,1)$
such that:%
\[
\sum_{\gamma\in\Gamma}\left(  1-\left\Vert \gamma\left(  0\right)  \right\Vert
_{\mathbb{C}^{n}}\right)  <\infty
\]
and the orbit $\Gamma(0)$ of $0$ spans $\mathbb{C}^{n}$. Then the automorphism
group of $F_{n,\Gamma(0)}^{\infty}=F_{n}^{\infty}/\mathcal{N}_{\Gamma(0)}$ is
the stabilizer subgroup $\underline{\Gamma}$ in $\operatorname*{SU}(n,1)$\ of
$\Gamma(0)$.\ Moreover, the action of $\Gamma$ on $F_{n,\Gamma(0)}^{\infty}$
is ergodic.
\end{theorem}

\begin{proof}
By Proposition (\ref{Blaschke}), the spectrum of $F_{n,\Gamma(0)}^{\infty}$ is
$\Gamma(0)$. By Theorem (\ref{Lift}), the group of automorphisms of
$F_{n,\Gamma(0)}^{\infty}$ is $\underline{\Gamma}$. To ease notation, let
$\tau$ be the group isomorphism from $\underline{\Gamma}$ onto
$\operatorname*{Aut}\left(  F_{n,\Delta}^{\infty}\right)  $ given by Theorem
(\ref{Lift}).

Now let $a\in F_{n,\Gamma(0)}^{\infty}$ such that for all $\gamma\in\Gamma$ we
have $\tau(\gamma)(a)=a$, so $\lambda\in\Gamma(0)\mapsto a(\lambda)$ is
constant, equal to $a(0)$. Hence $a-a(0)1$ is $0$ on $\Gamma(0)$ and thus
$\varphi-a(0)1\in\mathcal{N}_{\Gamma(0)}$, i.e. $a=\varphi+\mathcal{N}%
_{\Gamma(0)}=a(0)1$ is a scalar multiple of the identity. So $\Gamma$ acts ergodically.
\end{proof}

\bigskip We can deduce from Theorem (\ref{main}) the following simple
corollary, where the span condition is relaxed, but the relation between the
original group and the automorphism group may be less clear.

\begin{corollary}
Let $\Gamma$ be a a discrete subgroup of $\operatorname*{SU}(n,1)$ such that:%
\[
\sum_{\gamma\in\Gamma}\left(  1-\left\Vert \gamma\left(  0\right)  \right\Vert
_{\mathbb{C}^{n}}\right)  <\infty
\]
and the orbit $\Gamma(0)$ of $0$ spans $\mathbb{C}^{k}$ with $k\leq n$. Then
the automorphism group of $F_{n,\Gamma(0)}^{\infty}=F_{n}^{\infty}%
/\mathcal{N}_{\Gamma(0)}$ is isomorphic to the stabilizer subgroup
$\underline{\Gamma}$ in $\operatorname*{SU}(k,1)$\ of the image of $\Gamma(0)$
by a unitary $U$ of $\mathbb{C}^{n}$ such that $\operatorname*{Ad}U\left(
\Gamma\left(  0\right)  \right)  \subseteq\mathbb{C}^{k}\times\left\{
0\right\}  $.
\end{corollary}

\begin{proof}
Let $U$ be a unitary acting on $\mathbb{C}^{n}$ such that $\Delta
=U\Gamma(0)U^{\ast}\subseteq\mathbb{C}^{k}\times\left\{  0\right\}  $. Hence,
the span of $\Delta$ is $\mathbb{C}^{k}\times\left\{  0\right\}  $. Now, write
$\mathfrak{U}$ for the automorphism of $F_{n}^{\infty}$ such that
$\widehat{\mathfrak{U}}=\operatorname*{Ad}U$. Now $\mathfrak{U}$ maps
$\mathcal{N}_{\Gamma(0)}$ onto $\mathcal{N}_{\Delta}$ and thus defines an
isomorphism from $F_{n,\Gamma(0)}^{\infty}$ onto $F_{n,\Delta}^{\infty}$ still
denoted by $\mathfrak{U}$. It is now easy to see that the quotient map
$F_{n}^{\infty}\longrightarrow F_{n,\Delta}^{\infty}$ factors as
$F_{n}^{\infty}\longrightarrow F_{k}^{\infty}\longrightarrow F_{n,\Delta
}^{\infty}=F_{k,\Delta}^{\infty}$. Now by Theorem\ (\ref{Lift}), the group of
automorphisms of $F_{n,\Delta}^{\infty}$ is the stabilizer group
$\underline{\Gamma}$ of $\Delta$ in $\operatorname*{SU}\left(  k,1\right)  $.
Last, we note that $\Phi\in\operatorname*{Aut}\left(  F_{n,\Gamma(0)}^{\infty
}\right)  \mapsto\mathfrak{U}\circ\Phi\circ\mathfrak{U}^{-1}$ is a group
isomorphism onto $\operatorname*{Aut}\left(  F_{n,\Delta}^{\infty}\right)  $.
\end{proof}

\bigskip We conclude that, using methods similar to the proof of Theorem
(\ref{Lift}) we can also prove that:

\begin{corollary}
Let $\Delta_{1}$ and $\Delta_{2}$ be two subsets of $\mathbb{B}_{n}$. Then
there exists a completely isometric isomorphic from $F_{n,\Delta_{1}}^{\infty
}$ to $F_{n,\Delta_{2}}^{\infty}$ if and only if there exists an automorphism
$\gamma\in\operatorname*{SU}\left(  n,1\right)  $ of $\mathbb{B}_{n}$ such
that $\gamma\left(  \overline{\Delta_{1}}^{\Sigma}\right)  =\overline
{\Delta_{2}}^{\Sigma}$.
\end{corollary}

\bibliographystyle{amsplain}
\bibliography{thesis}

\end{document}